\documentclass[12pt]{amsart}
\usepackage{amsmath, amsthm, amscd, amsfonts}
\usepackage{color}

\newtheorem{theorem}{Theorem}[section]
\newtheorem{lemma}[theorem]{Lemma}
\newtheorem{proposition}[theorem]{Proposition}
\newtheorem{corollary}[theorem]{Corollary}
\theoremstyle{definition}
\newtheorem{definition}[theorem]{Definition}

\theoremstyle{remark}

\numberwithin{equation}{section}

\newfont{\kh}{msbm10}

\begin{document}
\title[Generic properties of module maps]
{Generic properties of module maps and characterizing inverse
limits of C*-algebras of compact operators}
\author{K. Sharifi}
\address{Kamran Sharifi,
Department of Mathematics, Shahrood University of Technology, P.
O. Box 3619995161-316, Shahrood, Iran
\newline School of Mathematics, Institute for Research in
 Fundamental Sciences (IPM), P.O. Box: 19395-5746, Tehran, Iran}
\email{sharifi.kamran@gmail.com and sharifi@shahroodut.ac.ir}

\thanks{This research was in part supported by a
grant from IPM (No. 90470018).}

\begin{abstract}

We study closedness of the range, adjointability and generalized
invertibility of modular operators between Hilbert modules over
locally C*-algebras of coefficients. Our investigations and the
recent results of M. Frank [Characterizing C*-algebras of compact
operators by generic categorical properties of Hilbert
C*-modules, {\it J. K-Theory} {\bf 2} (2008), 453-462] reveal a
number of equivalence properties of the category of Hilbert
modules over locally C*-algebras which characterize precisely the
inverse limit of C*-algebras of the C*-algebra of compact
operators.
\end{abstract}
\maketitle

\section{Introduction}
Locally C*-algebras are generalizations of C*-algebras. A locally
C*-algebra is a complete Hausdorff complex topological $*$-algebra
$ \mathcal{A}$, whose topology is determined by its continuous
C*-seminorms in the sense that the net $\{ a_i \}_{i \in I}$
converges to $0$ if and only if the net $\{p( a_i) \}_{i \in I}$
converges to $0$ for every continuous C*-seminorm $p$ on  $
\mathcal{A}$. Locally C*-algebras were first introduced by A.
Inoue \cite{Inoue1971} and studied more by N. C. Phillips and M.
Fragoulopoulou \cite{Fragoulopoulou2005, Phillips1988}. See also
the book of M. Joita \cite{JoitaBook1} and references therein.

Hilbert modules are essentially objects like Hilbert spaces by
allowing the inner product to take values in a (locally)
C*-algebra rather than the field of complex numbers. They play an
important role in the modern theory of operator algebras, in
noncommutative geometry and in quantum groups, see \cite{GVF}.

Throughout the present paper we refer to C*-subalgebras of the
C*-algebras of compact operators on Hilbert spaces as {\it
C*-algebras of compact operators}. Recall that a C*-algebra of
compact operators is a $c_{0}$-direct sum of elementary
C*-algebras $\mathcal{K}(H_{i})$ of all compact operators acting
on Hilbert spaces $H_{i}, \ i \in I$, cf.~\cite[Theorem
1.4.5]{ARV}.

Magajna and Schweizer, respectively, have shown that C*-algebras
of compact operators can be characterized by the property that
every closed (and coinciding with its biorthogonal complement,
respectively) submodule of every Hilbert C*-module over them is
automatically an orthogonal summand, cf. \cite{MAG, SCH}.
Together with results of Lj. Aramba\v{s}i\'c,  D. Baki\'c and B.
Gulja\v{s} \cite{ARA, B-G, GUL}, numerous generic properties of
the category of Hilbert C*-modules over C*-algebras which
characterize precisely the C*-algebras of compact operators have
been found by M. Frank and the author in \cite{FR1, F-S, FS2}.
The later work motivate us to study some properties of modular
operators, such as closedness of the range, adjointability, polar
decomposition and generalized invertibility of module maps
between Hilbert modules over locally C*-algebras of coefficients.
These help us to obtain a number of equivalence properties which
describe precisely the inverse limit of C*-algebras of compact
operators.

In the present paper we recall some definitions and simple facts
about Hilbert modules over locally C*-algebras and the module
maps between them. Then we study the closedness of the range and
adjointability of module maps, in fact we will prove that a
bounded module map between Hilbert modules over locally
C*-algebras is adjointable if and only if its graph is an
orthogonal summand (compare \cite{FR1}). A bounded adjointable
module map possesses a generalized inverse if and only if it has
a closed range. Finally, for a given locally C*-algebra $
\mathcal{A}$ we demonstrate that any bounded $\mathcal{A}$-module
map between arbitrary $\mathcal{A}$-modules possesses an adjoint
$\mathcal{A}$-module map, if and only if the images of all
bounded $ \mathcal{A}$-module maps with closed range between
arbitrary Hilbert $\mathcal{A}$-modules are orthogonal summands,
if and only if every bounded $\mathcal{A}$-module map between
arbitrary Hilbert $\mathcal{A}$-modules has polar decomposition,
if and only if every bounded $\mathcal{A}$-module map between
arbitrary Hilbert $\mathcal{A}$-modules has generalized inverse,
if and only if $\mathcal A$ is an inverse limit of C*-algebras of
compact operators.

\section{Preliminaries}
Suppose $\mathcal{A}$ is a locally C*-algebra and $S(
\mathcal{A})$ is the set of all continuous C*-seminorms on $
\mathcal{A}$. For every $p \in S( \mathcal{A})$, the quotient
$*$-algebra $ \mathcal{A}/N_{p}^{ \, \mathcal{A}} $ is denoted by
$ \mathcal{A}_p$\,, where $ N_{p}^{ \, \mathcal{A}} =\{ a \in
\mathcal{A}: p(a) = 0 \}$ is a C*-algebra in the C*-norm induced
by $p$. The canonical map from $ \mathcal{A}$ to  $\mathcal{A}_p$
is denoted by $ \pi_p^{ \, \mathcal{A}}$ and $a_p$ is reserved to
denote $\pi_p^{ \, \mathcal{A}}(a)$. For $p, q \in S(
\mathcal{A})$ with $p \geq q$, the surjective canonical map $
\pi_{pq}^{ \, \mathcal{A}} : \mathcal{A}_p \to \mathcal{A}_q$ is
defined by $ \pi_{pq}^{ \, \mathcal{A}} ( \pi_p^{ \,
\mathcal{A}}(a))= \pi_q^{ \, \mathcal{A}}(a)$ for all $a \in
\mathcal{A}$. Then $ \{ \mathcal{A}_p; \pi_{pq}^{ \, \mathcal{A}}
\}_{ p,\,q \in S( \mathcal{A}),\, p \geq q }$ is an inverse
system of C*-algebras and $\lim\limits_{\underset{p}{ \leftarrow
}} \mathcal{A}_{p}$ is a locally C*-algebra which can be
identified with $ \mathcal{A}$. We refer to the book
\cite{Fragoulopoulou2005} and papers \cite{Inoue1971,
Phillips1988} for  more information and useful examples.
 A morphism of locally C*-algebras
is a continuous $*$-morphism from a locally C*-algebra
$\mathcal{A}$ to another locally C*-algebra $ \mathcal{B}$. An
isomorphism of locally C*-algebras from $\mathcal{A}$ to $
\mathcal{B}$ is a bijective map $\Phi: \mathcal{A} \to
\mathcal{B}$ such that $\Phi$ and $\Phi^{-1}$ are morphisms of
locally C*-algebras.

A (right) {\it pre-Hilbert module} over a locally C*-algebra
algebra $\mathcal{A}$ is a right $\mathcal{A}$-module $E$,
compatible with the complex algebra structure, equipped with an
$\mathcal{A}$-valued inner product $\langle \cdot , \cdot \rangle
: E \times E \to \mathcal{A}\,, \ (x,y) \mapsto \langle x,y
\rangle$, which is $\mathcal A$-linear in the second variable $y$
and has the properties:
$$ \langle x,y \rangle=\langle y,x \rangle ^{*}, \ {\rm and} \
 \langle x,x \rangle \geq 0 \ \ {\rm with} \
   {\rm equality} \ {\rm if} \ {\rm and} \ {\rm only} \
   {\rm if} \ x=0.$$

A pre-Hilbert $\mathcal{A}$-module $E$ is a Hilbert
$\mathcal{A}$-module if $E$ is complete with respect to the
topology determined by the family of seminorms $ \{
\overline{p}_E \}_{p \in S( \mathcal{A})}$  where $\overline{p}_E
( \xi) = \sqrt{ p( \langle \xi, \xi \rangle)}$, $ \xi \in E$. If
$E$, $F$ are two Hilbert $ \mathcal{A}$-modules then the set of
all ordered pairs of elements $E \oplus F$ from $E$ and $F$ is a
Hilbert $\mathcal{A}$-module with respect to the $\mathcal
A$-valued inner product $\langle
(x_{1},y_{1}),(x_{2},y_{2})\rangle= \langle
x_{1},x_{2}\rangle_{E}+\langle y_{1},y_{2}\rangle _{F}$. It is
called the direct {\it orthogonal sum of $E$ and $F$}.

We say that a Hilbert $ \mathcal{A}$-submodule
$X$ of a Hilbert $ \mathcal{A}$-module $E$ is a topological summand if
$E$ can be decomposed into the
direct sum of the Banach $ \mathcal{A}$-submodule $X$ and of another Banach
$ \mathcal{A}$-submodule $Y$. The notation is $E=X \stackrel{.}{+} Y$.
If, moreover, the decomposition can be arranged as an orthogonal one (i.e. $X \perp Y$) then
the Hilbert $ \mathcal{A}$-submodule $X$ is an orthogonal summand of the
Hilbert $ \mathcal{A}$-module $E$. In this case, we write
$E=X \oplus Y$ and $Y=X^{ \perp}$.

Let $E$ be a Hilbert $\mathcal{A}$-module and $p\in
S(\mathcal{A})$, then $ N_{p}^{E}=\{\xi \in E; \
\overline{p}_{E}(\xi )=0\}$ is a closed submodule of $E$  and
$E_{p}=E / N_{p}^{E}$ is a Hilbert
$\mathcal{A}_{p}$-module with $(\xi
+ N_{p}^{E}{})\pi _{p}^{ \, \mathcal{A}}(a)=\xi a+N_{p}^{E}$ and $%
\left\langle \xi +N_{p}^{E}, \eta + N_{p}^{E} \right\rangle =\pi_{p}^{ \,
\mathcal{A}}(\left\langle \xi ,\eta \right\rangle ).$ The
canonical map from $E$ onto $E_{p}$ is denoted by $\sigma
_{p}^{E}$ and $\xi_p$ is reserved to denote $\sigma
_{p}^{E}(\xi)$. For $ p,q\in S( \mathcal{A})$ with $p\geq q$, the
surjective canonical map $ \sigma _{pq}^{E}
: E_{p} \to E_{q}$ is defined by $\sigma _{pq}^{E}(\sigma
_{p}^{E}(\xi ))=\sigma _{q}^{E}(\xi )$ for all $\xi \in E$. Then
$\{E_{p};\mathcal{A}_{p};\sigma _{pq}^{E},\pi _{pq}^{ \,
\mathcal{A}} \}_{p,\,q\in S(\mathcal{A}),\,p\geq q}$ is an
inverse system of Hilbert C*-modules in the following sense:
\begin{itemize}
\item $\sigma _{pq}^{E}(\xi _{p}a_{p})=\sigma _{pq}^{E}(\xi _{p})\pi
_{pq}^{ \, \mathcal{A}}(a_{p}),\ \xi _{p}\in E_{p}, \ a_{p}\in
\mathcal{A}_{p}, \ p,q\in S(\mathcal{A}), \ p\geq q$,

\item $\left\langle
\sigma _{pq}^{E}(\xi _{p}),\sigma _{pq}^{E}(\eta
_{p})\right\rangle =\pi _{pq}^{ \, \mathcal{A}}(\left\langle \xi
_{p},\eta _{p}\right\rangle ), \ \xi _{p},  \eta _{p} \in E_{p},
\ p,q\in S(\mathcal{A}), \ p\geq q$,

\item $\sigma _{qr}^{E} \circ \sigma _{pq}^{E}=\sigma
_{pr}^{E} $ if $p, q, r \in S(\mathcal{A})$, $p\geq q\geq r$, and

\item  $\sigma _{pp}(\xi _{p})=\xi _{p},\;\xi
_{p}\in E_{p}$, $p \in S(\mathcal{A})$.
\end{itemize}
In this case, $\lim\limits_{\underset{p}{\leftarrow }}E_{p}$ is a
Hilbert $\mathcal{A}$-module which can be identified with $E$.

Let $E$ and $F$ be Hilbert $ \mathcal{A}$-modules and $T : E \to
F$ be an $ \mathcal{A}$-module map. The module map $T$ is called
{\it bounded}  if for each $p \in S( \mathcal{A})$, there is $K_p
> 0$ such that $ \overline{p}_F(Tx) \leq K_p \,
\overline{p}_E(x)$  for all $x \in E$. The module map $T$ is
called {\it adjointable} if there exists an $ \mathcal{A}$-module
map $T^*: F \to E$ with the property $ \langle Tx,y \rangle
=\langle x,T^*y \rangle$ for all $x \in E$, $y \in F$. It is well
known that every adjointable $ \mathcal{A}$-module map is
bounded, cf. \cite[Lemma 2.2.3]{JoitaBook1}. The set
$\mathcal{L}_{\mathcal{A}}(E,F)$ of all bounded adjointable $
\mathcal{A}$-module maps from $E$ into $F$ becomes a locally
convex space with topology defined by the family of seminorms $
\{ \tilde{p}_{ \mathcal{L}_{\mathcal{A}}(E,F)} \}_{p \in S(
\mathcal{A})}$, in which, $ \tilde{p}_{
\mathcal{L}_{\mathcal{A}}(E,F)}(T)= \| (\pi_p^{ \,
\mathcal{A}})_{*}(T) \|_{\mathcal{L}_{\mathcal{A}_p}(E_p,F_p)} $
and $ (\pi_p^{ \, \mathcal{A}})_{*} :
\mathcal{L}_{\mathcal{A}}(E,F) \to
\mathcal{L}_{\mathcal{A}_p}(E_p,F_p)$ is defined by $(\pi_p^{ \,
\mathcal{A}})_{*}(T)( \xi + N_{p}^{E})=T \xi + N_{ p}^{F} $ for
all $T \in \mathcal{L}_{\mathcal{A}}(E,F)$, $ \xi \in E$. Suppose
$ p,q\in S(\mathcal{A}), \ p\geq q$ and $ ( \pi_{pq}^{ \,
\mathcal{A}})_{*}: \mathcal{L}_{ \mathcal{A}_p}(E_p,F_p) \to
\mathcal{L}_{ \mathcal{A}_q}(E_q,F_q)$ is defined by $(\pi_{pq}^{
\, \mathcal{A}})_{*}(T_p)( \sigma_{q}^{E}( \xi)) = \sigma
_{pq}^{F}(T_p( \sigma_{p}^{E}( \xi))$. Then $ \{ \mathcal{L}_{
\mathcal{A}_p}(E_p,F_p); \ (\pi_{pq}^{ \, \mathcal{A}})_{*} \}_{
p,q\in S(\mathcal{A}), \, p\geq q}$ is an inverse system of Banach
spaces and $\lim\limits_{\underset{p}{ \leftarrow }}
\mathcal{L}_{ \mathcal{A}_p}(E_p,F_p)$ can be identified by
$\mathcal{L}_{\mathcal{A}}(E,F)$. In particular, topologizing,
$\mathcal{L}_{\mathcal{A}}(E,E)$  becomes a locally C*-algebra
which is abbreviated by $\mathcal{L}_{\mathcal{A}}(E)$. Proofs of
the above facts can be founded in Sections 2.1 and 2.2 of the book
\cite{JoitaBook1}. Hilbert modules over locally C*-algebras have
been studied systematically in the book \cite{JoitaBook1} and the
papers \cite{Joita1, Joita2, Joita3, Phillips1988}.

We use the notations $Ker(\cdot)$ and $Ran(\cdot)$ for kernel and
range of module maps, respectively. A bounded $
\mathcal{A}$-module map $P:E \to E$ is said to be idempotent if
$P^2=P$. If, in addition, $P$ is adjointable and $P^*=P$ then $P$
is said to be projection. It is known that a Hilbert $
\mathcal{A}$-submodule $X$ of a Hilbert $ \mathcal{A}$-module $E$
is an orthogonal summand (a topological summand, respectively) if
and only if there exists a projection (an idempotent,
respectively) on $E$ whose range is $X$.
\begin{lemma}\label{proj1}Suppose $P:E \to E$ is a bounded $ \mathcal{A}$-module map. Then
$P$ is an idempotent if and only if $ (\pi_{p}^{ \,
\mathcal{A}})_{*}(P): E_p \to E_p$, $(\pi_p^{ \,
\mathcal{A}})_{*}(P)( \xi + N_{p}^{E})=P \xi + N_{ p}^{E} $ is an
idempotent for each $p \in S( \mathcal{A})$. In particular, $P$ is
a projection in $\mathcal{L}_{\mathcal{A}}(E)$ if and only if $
(\pi_{p}^{ \, \mathcal{A}})_{*}(P)$ is a projection in
$\mathcal{L}_{\mathcal{A}_p}(E_p)$ for each $p \in S(
\mathcal{A})$.
\end{lemma}

\begin{proof}Suppose $P$ is an idempotent and $p \in S(
\mathcal{A})$. Then $ (\pi_{p}^{ \, \mathcal{A}})_{*}(P): E_p \to
E_p$ is a bounded $\mathcal{A}_p$-module map and for each
$x_p, y_p \in E_p$ we have $$ ((\pi_{p}^{ \,
\mathcal{A}})_{*}(P))^2 \, x_p= P^2x + N_p^{E}= Px + N_p^{E}=
(\pi_{p}^{ \, \mathcal{A}})_{*}(P) \, x_p,$$ that is, $
(\pi_{p}^{ \, \mathcal{A}})_{*}(P)$ is an idempotent.

Conversely, suppose $ (\pi_{p}^{ \, \mathcal{A}})_{*}(P) : E_p \to E_p$
is an idempotent. We obtain
\begin{eqnarray*}
\pi_{p}^{ \, \mathcal{A}}( \langle P^{2}x,y  \rangle - \langle Px,
y \rangle  ) &=& \langle  (\pi_{p}^{ \, \mathcal{A}})_{*}(P^2) \,
x_p, y_p \rangle - \langle (\pi_{p}^{ \, \mathcal{A}})_{*}(P) \,
x_p, y_p  \rangle \\ &=& \langle  ((\pi_{p}^{ \,
\mathcal{A}})_{*}(P))^2 \, x_p, y_p \rangle - \langle (\pi_{p}^{
\, \mathcal{A}})_{*}(P) \, x_p, y_p  \rangle = 0
\end{eqnarray*}
for  all $p \in S( \mathcal{A})$ and $x, y \in E$. We therefore
have  $\langle P^{2}x,y  \rangle = \langle Px, y \rangle,$ i.e.,
$P$ is an idempotent. A similar argument shows that $P$ is
selfadjoint if and only if $ (\pi_{p}^{ \, \mathcal{A}})_{*}(P)$
is. This proves the second statement.
\end{proof}

\begin{corollary}\label{proj0-0}Suppose $F$ and $E$ are Hilbert
$ \mathcal{A}$-modules which are identified with
$\lim\limits_{\underset{p}{ \leftarrow }} F_{p}$ and
$\lim\limits_{\underset{p}{ \leftarrow }} E_{p}$, respectively. If
$E$ is a $ \mathcal{A}$-submodule of $F$, $E$ is topologically
(orthogonally) complemented if and only if $E_p$ is topologically
(orthogonally) complemented for each $p \in S( \mathcal{A})$.
\end{corollary}

\begin{lemma}\label{proj2}Let $T$ be a module map in
$ \mathcal{L}_{\mathcal{A}}(E,F)$ which can be identified by $(
T_p )_p$ in $\lim\limits_{\underset{p}{ \leftarrow }}
\mathcal{L}_{ \mathcal{A}_p}(E_p,F_p)$. Then $T$ has closed range
if and only if $ (\pi_{p}^{ \, \mathcal{A}})_{*}(T)$ has closed
range for each $p \in S( \mathcal{A})$.
\end{lemma}

\begin{proof} For $P \in \mathcal{L}_{\mathcal{A}}(F)$,
$Ran(P)=Ran(T)$ if and only if $Ran((\pi_{p}^{ \,
\mathcal{A}})_{*}(P))=Ran((\pi_{p}^{ \, \mathcal{A}})_{*}(T))$
for each $p \in S( \mathcal{A})$. The result follows the above
fact, Lemma \ref{proj1} and \cite[Theorem 2.2]{Joita1}.
\end{proof}

Closed submodules of Hilbert modules need not be orthogonally
complemented at all, but Theorem 2.2 of \cite{Joita1}, which is
an extension of \cite[Theorem 3.2]{LAN}, states under which
conditions closed submodules may be orthogonally complemented.
For the special choice of modular operator $T \in
\mathcal{L}_{\mathcal{A}}(E,F)$ with closed range one has:

\begin{itemize}
\item $Ker(T)$ is orthogonally complemented in $E$, with
      complement $Ran(T^*)$,
\item $Ran(T)$ is orthogonally complemented in $F$, with
      complement $Ker(T^*)$,
\item the map $T^* \in \mathcal{L}_{\mathcal{A}}(F,E)$ has a closed range, too.
\end{itemize}

An $ \mathcal{A}$-module map $U \in
\mathcal{L}_{\mathcal{A}}(E,F)$ is said to be unitary if $U^*U
=1_E$ and $UU^* =1_F$. If there exists a unitary element of
$\mathcal{L}_{\mathcal{A}}(E,F)$ then we say that $E$ and $F$ are
unitarily equivalent Hilbert $\mathcal{A}$-modules. Two Hilbert
$\mathcal{A}$-modules $E$ and $F$ are isomorphic if and only if
there is a unitary operator from $E$ to $F$, cf. \cite[Corollary
2.5.4]{JoitaBook1}.

\begin{lemma}  \label{JoitaRemark2.5.2}{\rm (cf. \cite[Remark 2.5.2]{JoitaBook1})}
Suppose $U \in \mathcal{L}_{\mathcal{A}}(E,F)$. Then $U$ is
unitary if and only if $ (\pi_{p}^{ \, \mathcal{A}})_{*}(U) : E_p
\to F_p$ is a unitary operator for all $p \in S( \mathcal{A})$.
\end{lemma}

Let $E, F$ be $ \mathcal{A}$-modules and $T: E \to F$ be an $
\mathcal{A}$-module map then $ \mathcal{A}$-submodule
$G(T)=\{(x,Tx): \ x \in E \}$ is called the {\it graph} of $T$.
If $T$ is bounded $ \mathcal{A}$-module map then $G(T)$ is a
closed $ \mathcal{A}$-submodule of the Hilbert $
\mathcal{A}$-module $E \oplus F$. It is well known that a bounded
module map between Hilbert C*-modules is adjointable if and only
if its graph is an orthogonal summand, see e.g., \cite{FR1}. The
problem are restudied in the case of unbounded module maps
between Hilbert C*-modules in \cite{F-S}. In this section we
study adjointability of bounded $ \mathcal{A}$-module maps
between Hilbert modules over locally C*-algebras. In fact, we
show that the Hilbert $ \mathcal{A}_p$-modules $G(T)_p$ and $G(
(\pi_{p}^{ \, \mathcal{A}})_{*}(T) )$ are isomorphic and then we
lift Corollary 2.4 of \cite{FR1} to the case of Hilbert modules
over locally C*-algebras.

\begin{lemma}\label{proj3}Suppose $T: E \to F$ is bounded $
\mathcal{A}$-module map then the Hilbert $
\mathcal{A}_p$-modules $G(T)_p$ and $G( (\pi_{p}^{ \,
\mathcal{A}})_{*}(T) )$ are isomorphic for every $p \in S(
\mathcal{A})$.
\end{lemma}
\begin{proof} Suppose $p \in S(
\mathcal{A})$ and $T_p=(\pi_{p}^{ \, \mathcal{A}})_{*}(T)$. We
define the $ \mathcal{A}$-module maps $$U_p: G(T)_p \to
G(T_p), ~~U_p((x,Tx)+N_p^{G(T)})=(x_p, T_{p\,} x_p) $$ and $$W_p
: G(T_p) \to G(T)_p,~~W_p(x_p, T_{p\,} x_p)=
((x,Tx)+N_p^{G(T)}),$$ we obtain
$$ \langle U_p((x,Tx)+N_p^{G(T)}), (y_p, T_{p\,} y_p)  \rangle =
\langle (x,Tx)+N_p^{G(T)}, W_p(y_p, T_{p\,} y_p) \rangle$$ for
all $x \in E$, $y_p \in E_p$. That is, $U_p$ is adjointable and
$U_p^{*}=W_p$. We also have $U_p U_p^{*} =  U_p^{*} U_p=1$ on
$G(T)_p$, i.e., $G(T)_p$ and $G(T_p)$ are unitarily equivalent.
The result now follows from Corollary 2.5.4 of \cite{JoitaBook1}.
\end{proof}
\begin{lemma}\label{proj4}A bounded $ \mathcal{A}$-module map $T: E \to F$
is adjointable if and only if $ (\pi_{p}^{ \,
\mathcal{A}})_{*}(T): E_{p} \to F_{p}$ is adjointable for each $p
\in S( \mathcal{A})$. In this situation, the adjoint of $(\pi_{p}^{ \,
\mathcal{A}})_{*}(T)$ is $(\pi_{p}^{ \, \mathcal{A}})_{*}(T^*).$
\end{lemma}
\begin{proof} Suppose $T: E \to F$ is adjointable then $ \pi_p^{
\, \mathcal{A}}( \langle Tx,y  \rangle)=\pi_p^{ \, \mathcal{A}}(
\langle x,T^*y  \rangle)$ for all $x \in E$, $y \in F$ and $p \in
S( \mathcal{A})$. We therefore have $ \langle (\pi_{p}^{ \,
\mathcal{A}})_{*}(T) \, x_p, y_p   \rangle = \langle  x_p,
(\pi_{p}^{ \, \mathcal{A}})_{*}(T^*) \, y_p \rangle $, for all
$x_p \in E_p$, $y_p \in F_p$ and $p \in S( \mathcal{A})$, i.e.,
$(\pi_{p}^{ \, \mathcal{A}})_{*}(T)$ is adjointable and its
adjoint is $(\pi_{p}^{ \, \mathcal{A}})_{*}(T^*)$.

Conversely, suppose $T_p=(\pi_{p}^{ \, \mathcal{A}})_{*}(T):
E_{p} \to F_{p}$ is adjointable for all $p \in S( \mathcal{A})$.
Suppose $S: F \to E$ is defined by $Sy=( T_p^{ \,*} \,
y_p)_{p}$\,,  $y=(y_p)_p \in F= \lim\limits_{\underset{p}{
\leftarrow }} F_{p}$. Then $S$ is well defined, since
$$ \sigma _{pq}^{E}( \, T_{p}^{*} \, y_p)= (\pi _{pq}^{\,
\mathcal{A}} )_{*}( T_{p}^{*})(\sigma _{pq}^{F}(y_p))=T_{q}^{*}
\,y_q $$ for all $p,q \in S( \mathcal{A})$ with $p\geq q$.
Furthermore, we have
$$ \pi_p^{ \, \mathcal{A}}( \langle Tx,y  \rangle)= \langle T_p \, x_p, y_p \rangle
= \langle x_p, T_p^{ \, *} \, y_p \rangle  =\pi_p^{ \,
\mathcal{A}}( \langle x,Sy  \rangle) $$ for all $x=(x_p)_p \in E$,
$y=(y_p)_p \in F$ and $p \in S(\mathcal{A})$. Therefore
$ \langle Tx,y  \rangle = \langle x,Sy
  \rangle $ , i.e., $T$ is adjointable and $S=T^*$.
\end{proof}

\begin{proposition}\label{proj5}A bounded $ \mathcal{A}$-module map $T: E \to F$
possesses an adjoint map $T^*: F \to E$ if and only if the graph
of $T$ is an orthogonal summand of the Hilbert $
\mathcal{A}$-module $ E \oplus F$.
\end{proposition}

\begin{proof} Using Lemmas  \ref{proj3}, \ref{proj4} and  \cite[Corollary 2.4]{FR1},
we conclude that $T$ is adjointable if and only if every $
(\pi_{p}^{ \, \mathcal{A}})_{*}(T)$ is adjointable, if and only
if every $G( (\pi_{p}^{ \, \mathcal{A}})_{*}(T) )$ is an
orthogonal summand, if and only if $G(T)_p$ is an orthogonal
summand.

According Corollary \ref{proj0-0} and the fact that $G(T)=
\lim\limits_{\underset{p}{ \leftarrow }} G(T)_{p}$\,, the $
\mathcal{A}_p$-submodule $G(T)_p$ is an orthogonal summand in $(E
\oplus F)_p$ if and only if the $ \mathcal{A}$-submodule  $G(T)$
is an orthogonal summand in $E \oplus F$, which completes the
proof.

\end{proof}

\section{Polar decomposition and generalized inverses}
The polar decomposition is a useful tool that represents an
operator as a product of a partial isometry and a positive
element. It is well known that every bounded operator on Hilbert
spaces has polar decomposition. In general bounded adjointable
$\mathcal{A}$-module maps between Hilbert $\mathcal{A}$-modules
do not have polar composition, but M. Joita has given a necessary
and sufficient condition for bounded adjointable module maps to admit polar
decomposition. She has proved that a bounded adjointable operator
$T$ has polar decomposition if and only if $\overline{Ran(T)}$ and
$\overline{Ran(|T|)}$ are orthogonal direct summands. The reader
is encouraged to see \cite[Theorem 2.8, Proposition 2.10]{Joita1}
and \cite[Section 3.3]{JoitaBook1} for more information and the
proof of this fact.  See also Theorem 15.3.7 of \cite{WEG}.

\begin{definition} An adjointable module map $T: E \to F$ has a polar
decomposition if there is a partial isometry $V:E \to F$ such
that $T = V |T|$,  and $Ker (V) = Ker (T)$, $Ran(V) = \overline{Ran(T)}$,
$Ker (V^*)  = ker (T^*)$ and $Ran(V^*)= \overline{Ran(|T|)}$.
\end{definition}

\begin{proposition}\label{proj6} A bounded adjointable $ \mathcal{A}$-module map
 $T: E \to F$
 has a polar decomposition if and only if $(\pi_{p}^{ \, \mathcal{A}})_{*}(T) $
has a polar decomposition for each $p \in S( \mathcal{A})$. In
this situation, $T = V |T|$ if and only if $(\pi_{p}^{ \,
\mathcal{A}})_{*}(T) = (\pi_{p}^{ \, \mathcal{A}})_{*} (V) |
(\pi_{p}^{ \, \mathcal{A}})_{*}(T)| $ for each $p \in S(
\mathcal{A})$.
\end{proposition}

\begin{definition}Let $T \in \mathcal{L}_{\mathcal{A}}(E,F)$, then a bounded
adjointable operator $T^{ \dag} \in
\mathcal{L}_{\mathcal{A}}(F,E)$ is called the {\it generalized
inverse} of $T$ if
\begin{equation} \label{MP inverse}
T \, T^{ \dag}T=T, \ T^{ \dag}T \, T^{ \dag}= T^{ \dag}, \ (T \,
T^{ \dag})^*=T \, T^{ \dag} \ {\rm and} \ ( T^{ \dag} T)^*= T^{
\dag} T.
\end{equation}
\end{definition}

 The notation $T^{ \dag}$ is reserved to denote the
generalized inverse of $T$. These properties imply that $T^{
\dag}$ is unique and $ T^{ \dag} T$ and $ T \, T^{ \dag} $ are
orthogonal projections. Moreover, $Ran( T^{ \dag} )=Ran( T^{
\dag}  T)$, $Ran( T )=Ran( T \, T^{ \dag})$,  $Ker(T)=Ker( T^{
\dag} T)$ and $Ker(T^{ \dag})=Ker( T \, T^{ \dag} )$ which lead
us to $ E= Ker( T^{ \dag} T) \oplus Ran( T^{ \dag} T)= Ker(T)
\oplus Ran( T^{ \dag} )$ and $F= Ker(T^{ \dag}) \oplus Ran(T).$

 Xu and Sheng in \cite{Xu/Sheng} have shown that a bounded
adjointable operator between two Hilbert C*-modules admits a
bounded generalized inverse if and only if the operator has
closed range. The reader should be aware of the fact that a
bounded adjointable operator may admit an unbounded operator as
its generalized, see \cite{FS2, SHA/PARTIAL, SHA/Groetsch} for
more detailed information.

\begin{lemma}\label{proj6-7}
Let $T \in \mathcal{L}_{\mathcal{A}}(E,F) $, then $T$ has closed
range if and only if $Ker(T)$ is orthogonally complemented in $E$
and $T$ is bounded below on $Ker(T)^{\perp}$, i.e. for each $p
\in S( \mathcal{A})$ there exist $c_p
> 0$ such that $\overline{p}_F(Tx) \geq c_p \,
\overline{p}_E(x)$,  for all $x \in Ker(T)^{\perp}$. In this
case, $ (T_{| Ker(T)^{ \perp}})^{-1}$ is a bounded module map on
$Ran(T)$.
\end{lemma}
\begin{proof} Let first $Ran(T)$ be closed then $Ker(T)$ is orthogonally complemented in $E$.
Identifying $T$ with $(T_p) \in \lim\limits_{\underset{p}{
\leftarrow }} \mathcal{L}_{ \mathcal{A}_p}(E_p,F_p)$, then for
every $ p \in S( \mathcal{A})$ the range of $T_p$ is closed. According to
\cite[Proposition 1.3]{F-S}, $Ker(T_p)$ is orthogonally complemented and there exists
$c_p > 0$ such that $ \| T_p \, x_p \| \geq c_p \| x_p \| $ for
all $x \in Ker(T_p)^{\perp}$.
The latter inequality implies that
\begin{eqnarray*}
\overline{p}_F(Tx)^2 &=& p( \langle Tx,Tx \rangle) =
\| \pi_{p}^{ \mathcal{A}}( \langle Tx,Tx \rangle)\\
&=&\| \langle \sigma_{p}^{F}(Tx),\sigma_{p}^{F}(Tx) \rangle \|\\
&=&\| \langle T_p(\sigma_{p}^{E}(x)), T_p(\sigma_{p}^{E}(x)) \rangle \|\\
& \geq &c_{p}^2 \, \| \langle \sigma_{p}^{E}(x), \sigma_{p}^{E}(x) \rangle \|\\
& = &c_{p}^2 \, \| \pi_{p}^{ \mathcal{A}}( \langle x,x \rangle)
\|= c_{p}^2 \, \overline{p}_{E}(x)^2.
\end{eqnarray*}
Consequently, for each $p \in S( \mathcal{A})$ there exists $c_p >
0$ such that $\overline{p}_F(Tx) \geq c_p \, \overline{p}_E(x)$,
for all $x \in Ker(T)^{\perp}$.

The converse can be proved by a similar manner as the proof of
\cite[Proposition 1.3]{F-S}, and so we omitted it. The second
assertion follows from the first assertion.
\end{proof}

\begin{proposition}\label{proj7}Suppose $T \in \mathcal{L}_{\mathcal{A}}(E,F)$. The
operator $T$ has a generalized inverse if and only if $T$ has a
closed range.
\end{proposition}

\begin{proof} Suppose $T$ has a generalized inverse, then $ T \, T^{ \dag}
$ is an orthogonal projection which implies the closedness of $Ran(
T )=Ran( T \, T^{ \dag})$.

Conversely, suppose $Ran(T)$ is closed then $E=Ker(T) \oplus Ran(T^*)$ and
$F=Ker(T^*) \oplus Ran(T)$ by \cite[Theorem 2.2]{Joita1}. According to
Lemma \ref{proj6-7}, the module maps $T_{| Ker(T)^{ \perp}}$ and
$T^{*}_{~~| Ker(T^*)^{ \perp}}$ are invertible on $Ran(T)$ and $Ran(T^*)$, respectively,
which allow us to define $ \mathcal{A}$-module map $T^{ \dag} : F \to E$ and
$T^{ \dag \, *} : E \to F$ by

\[ \qquad T^{ \dag} x \, =~
\begin{cases}
(T_{| Ker(T)^{ \perp}})^{-1}x ~ & \text{ if   $x \in Ran(T)$}\\
0 ~ & \text{ if $x \in Ker(T^*), $}\
\end{cases}
\] \\
\[ \qquad T^{ \dag \, *}x =
\begin{cases}
(T^{*}_{~~| Ker(T^*)^{ \perp}})^{-1}x ~ & \text{ if   $x \in Ran(T^*)$}\\
0 ~ & \text{ if $x \in Ker(T) $}.\
\end{cases}
\]
Using the orthogonal direct sum decompositions, the module maps
$T^{ \dag}$ and $T^{ \dag \, *}$ satisfy $ \langle T^{ \dag} x,y
\rangle = \langle x, T^{ \dag \, ^*}y \rangle $ for all $x \in F$
and $y \in E$, which implies that $T^{ \dag} \in
\mathcal{L}_{\mathcal{A}}(F,E)$. Moreover, $T$ and $T^{ \dag}$
satisfy (\ref{MP inverse}), i.e., $T^{ \dag}$ is the generalized
inverse of $T$.
\end{proof}

\begin{corollary}\label{proj9}Suppose $T \in \mathcal{L}_{\mathcal{A}}(E,F)$.
The module map $T$ has generalized inverse if and only if $
(\pi_{p}^{ \, \mathcal{A}})_{*}(T)$ has generalized inverse for
each $p \in S( \mathcal{A})$. In this case, $ ( \, (\pi_{p}^{ \,
\mathcal{A}})_{*}(T))^{ \dag} = (\pi_{p}^{ \,
\mathcal{A}})_{*}(T^{ \dag})$ for each $p \in S( \mathcal{A})$.
\end{corollary}

The above result follows from the previous proposition, Lemma
\ref{proj2} and \cite[Theorem 2.2]{Xu/Sheng}. Let $ \mathcal{A}$
be a locally C*-algebra and $ \mathfrak{a} \in \mathcal{A}$. An
element $ \mathfrak{a}^{ \dag} \in \mathcal{A}$ is called the
generalized inverse of $\mathfrak{a}$ if $\mathfrak{a}$ and
$\mathfrak{a}^{ \dag}$ satisfy (\ref{MP inverse}). Generalized
inverses in C*-algebras have been investigated by R. Harte and M.
Mbekhta \cite{Harte/Mbekhta}. The main result of their paper now
reads as follows:

\begin{corollary}\label{proj8}Suppose $ \mathcal{A}$ is a unital
locally C*-algebra and $ \mathfrak{a} \in \mathcal{A}$. Then
$\mathfrak{a}$ has a generalized inverse if and only if $
\mathfrak{a} \, \mathcal{A}$ is a closed right ideal in $
\mathcal{A}$.
\end{corollary}

Since every locally C*-algebra is a right $ \mathcal{A}$-module
on its own, the fact directly follows from Proposition
\ref{proj7}.

\section{Inverse limits of C*-algebras of compact
operators}

We closed the paper with characterizing the inverse limit of
C*-algebras of compact operators via the generic properties of
module maps. To deduce the following theorem just one needs to
use \cite[Theorem 2.6]{FR1}, Corollary \ref{proj0-0}, Lemmas
\ref{proj2}, \ref{proj4}, \ref{proj9} and
Proposition \ref{proj6}.\\

\begin{theorem}\label{proj0}
Let $\mathcal{A}$ be a locally C*-algebra. The following
conditions are equivalent:
\newcounter{cou001}
\begin{list}{(\roman{cou001})}{\usecounter{cou001}}
\item $\mathcal{A}$ is an inverse limit of C*-algebras of compact operators.

\item For every Hilbert $\mathcal{A}$-module $E$ every
Hilbert $\mathcal{A}$-submodule $F \subseteq E$ is automatically
orthogonally complemented, i.e. $F$ is an orthogonal summand.

\item For every Hilbert $\mathcal{A}$-module $E$ Hilbert
$\mathcal{A}$-submodule $F \subseteq E$ that coincides with its
biorthogonal complement $F^{\perp \perp} \subseteq E$ is
automatically orthogonally complemented in $E$.

\item For every pair of Hilbert $\mathcal{A}$-modules
$E, F$, every bounded $\mathcal{A}$-module map $T:E \to F$
possesses an adjoint bounded $\mathcal{A}$-module map $T^*:F \to
E$.

\item The kernels of all bounded $ \mathcal{A}$-module
maps between arbitrary Hilbert $ \mathcal{A}$-modules are
orthogonal summands.

\item The image of all bounded $ \mathcal{A}$-module maps
with norm closed range between arbitrary Hilbert
$\mathcal{A}$-modules are orthogonal summands.

\item For every pair of Hilbert $ \mathcal{A}$-modules
$E, F$, every bounded $\mathcal{A}$-module map $T : E \to F$ has
polar decomposition, i.e. there exists a unique partial isometry
$V$ with initial set $ \overline{Ran(|T|)}$ and the final set $
\overline{Ran(T)}$ such that $T=V|T|$.

\item For every pair of Hilbert $ \mathcal{A}$-modules
$E, F$, every bounded $\mathcal{A}$-module map $T: E \to F$ has
generalized inverse.
\item For every Hilbert $ \mathcal{A}$-module $E$
every Hilbert $ \mathcal{A}$-submodule is automatically
topologically complemented there, i.e. it is a topological
direct summand.
\end{list}
\end{theorem}

{\bf Acknowledgement}: The author would like to thank professor
M. Joita who sent the author some copies of her recent
publications. The author is also grateful to the referee for
his/her careful reading and his/her useful comments.


\begin{thebibliography}{99}


\bibitem {ARA} Lj. Aramba\v{s}i\'c, Another characterization of Hilbert C*-modules
over compact operators, {\it J. Math. Anal. Appl.}  {\bf 344}
(2008), no. 2, 735-740.

\bibitem {ARV} W. Arveson, {\it An Invitation to C*-algebras},
  Springer, New York, 1976.


\bibitem {B-G} D. Baki\'c and B. Gulja\v{s}, Hilbert C*-modules
over C*-algebras of compact operators, {\it Acta Sci. Math.}
(Szeged) {\bf 68} (2002), no. 1-2, 249-269.

\bibitem {FR2} M. Frank, Geometrical aspects of Hilbert C*-modules,
{\it Positivity} {\bf 3} (1999), 215-243.

\bibitem {FR1} M. Frank, Characterizing C*-algebras of compact
operators by generic categorical properties of Hilbert
C*-modules, {\it J. K-Theory} {\bf 2} (2008), 453-462.
\bibitem {F-S} M. Frank and K. Sharifi, Adjointability of densely
  defined closed operators and the Magajna-Schweizer theorem,
  {\it J. Operator Theory}  {\bf 63} (2010),  271-282.
\bibitem {FS2} M. Frank and K. Sharifi, Generalized inverses and polar
decomposition of unbounded regular operators on Hilbert
C*-modules, {\it J. Operator Theory} {\bf 64} (2010), 377-386.
\bibitem {Fragoulopoulou2005} M. Fragoulopoulou, {\it Topological algebras with involution},
 North Holland, Amsterdam, 2005.
\bibitem {GUL} B. Gulja\v{s}, Unbounded operators on Hilbert
C*-modules over C*-algebras of compact operators, {\it J. Operato
Theory}   {\bf 59}  (2008),  no. 1, 179-192.
\bibitem {GVF} J. M. Gracia-Bond\'{\i}a and J. C. V\'arilly and H. Figueroa,
{\it Elements of non-commutative geometry}, Birkh\"auser, 2000.
\bibitem {Harte/Mbekhta} R. Harte and M. Mbekhta, Generalized inverses in C*-algebras
II, { \it Studia Math.}, {\bf 103} (1992), 71-77.
\bibitem {H-A} S. Hejazian and T. Aghasizadeh, Equivalence classes of linear
mappings on $ \mathcal{B(M)}$, to appear in {\it Bull. Malays.
Math. Sci. Soc.} 
\bibitem{Inoue1971} A. Inoue, Locally C*-algebras,
{\it Mem. Faculty Sci. Kyushu Univ. Ser. A} {\bf 25}
(1971),197-235.
\bibitem{JoitaBook1}M. Joita, { \it Hilbert modules over locally
C*-algebras,} University of Bucharest Press, 2006.
\bibitem{Joita1}M. Joita, Projections on Hilbert modules over locally C*-algebras,
{ \it Math. Reports,} { \bf 54} (2002), 373-378.
\bibitem{Joita2}M. Joita, Multipliers of locally C*-algebras,
{ \it An. Univ. Bucuresti, Mat.} { \bf 48} (1999), 17-24.
\bibitem{Joita3}M. Joita, On Hilbert modules over locally C*-algebras,
{ \it An. Univ. Bucuresti, Mat.} { \bf 49} (2000), 41-51.
\bibitem{Joita4}M. Joita, On Hilbert modules over locally C*-algebras II,
{ \it Period. Math. Hungar.}, { \bf 51} (2005), 27-36.


\bibitem {LAN} E. C. Lance, {\it Hilbert C*-Modules}, LMS
Lecture Note Series 210, Cambridge Univ. Press, 1995.

\bibitem {MAG} M. Magajna, Hilbert C*-modules in which all
closed submodules are complemented, {\it Proc. Amer. Math. Soc.}
{\bf 125(3)} (1997), 849-852.

\bibitem {Phillips1988} N. C. Phillips, Inverse limits of C*-algebras, {\it J. Operator
Theory}  {\bf 19} (1988) 159-195.

\bibitem {SCH} J. Schweizer, A description of Hilbert C*-modules in which
all closed submodules are orthogonally closed, {\it Proc. Amer.
Math. Soc.} {\bf 127} (1999), 2123-2125.
\bibitem {SHA/PARTIAL} K. Sharifi, Descriptions of partial isometries on Hilbert
C*-modules, {\it Linear Algebra Appl.} {\bf  431} (2009), 883-887.

\bibitem {SHA/Groetsch} K. Sharifi, Groetsch's representation of  Moore-Penrose inverses and
ill-posed problems in Hilbert C*-modules, {\it J. Math. Anal.
Appl.} {\bf 365} (2010), 646-652.

\bibitem {WEG} N. E. Wegge-Olsen, {\it K-theory and $C^*$-algebras:
 a Friendly Approach}, Oxford University Press, Oxford, England, 1993.

\bibitem{Xu/Sheng} Q. Xu and L. Sheng, Positive semi-definite
matrices of adjointable operators on Hilbert C*-modules, {\it
Linear Algebra Appl.} {\bf428} (2008), 992-1000.

\end{thebibliography}
\end{document}